\documentclass[12pt]{amsart}
\usepackage{amsmath,amsthm,amssymb}
 \usepackage[all]{xy}

\def\pmatrix{\left(\begin{matrix}}
\def\endpmatrix{\end{matrix}\right)}
\def\s{\sigma}

\def\<{\langle}
\def\>{\rangle}

\def\GL{\operatorname{GL}}

\def\diag{\operatorname{diag}}
\def\Thr{\operatorname{Th}_r}
\def\Thtwo{\operatorname{Th}_2}
\def\IM{\operatorname{Im}}

\def\3{^{(3)}}
\def\4{^{(4)}}
\def\5{^{(5)}}

\def\Th2{\operatorname{Th}_2}

\let\th\theta
\def\({\left(}
\def\){\right)}

\def\P{{\mathbb P}}
\def\F{{\mathbb F}}
\def\Z{{\mathbb Z}}
\def\C{{\mathbb C}}

\def\diag{\operatorname{diag}}

\def\t{\theta}
\def\s{\sigma}

\def\H{{\mathcal H}}

\def\tch#1#2{{\left[\begin{matrix}{#1}\\ {#2}\end{matrix}\right]}}
\def\tt#1#2{{\t\tch#1#2}}

\theoremstyle{plain}
\newtheorem{thm}{Theorem}

\newtheorem{prop}[thm]{Proposition}

\newtheorem{cor}[thm]{Corollary}

\theoremstyle{definition}

\title{ On the
image of code polynomials under theta map}
\author{Manabu Oura, Riccardo Salvati Manni}
\address{Department of Mathematics,
Faculty of Science,
Kochi University,
Akebono-cho 2-5-1,
Kochi 780-8520,
Japan}
\email{oura@math.kochi-u.ac.jp}
\address{Dipartimento di Matematica, Universit\`a ``La Sapienza'',
Piazzale A. Moro 2, Roma, I 00185, Italy}
\email{salvati@mat.uniroma1.it}
\begin{document}
\begin{abstract}
The theta map sends code polynomials into the ring of Siegel modular
forms of even weights.
Explicit description of the image is known for $g\leq 3$
and the surjectivity of the theta map follows.
Instead it is known that this map is not surjective for $g\geq 5$.
In this paper we discuss the possibility of an embedding between the 
associated projective varieties. We prove that this is not possible 
for $g\geq 4$ and consequently we get  the non surjectivity of the graded rings
for the remaining case $g=4$.  
\end{abstract}
\date{\today}
 
\maketitle

\section{Introduction}
One of the main theme in the theory of modular forms is 
to determine explicitly the structure of 
the graded rings of modular forms  in terms of generators and
relations. 
This study started systematically with Igusa about 50 years ago.\smallskip

In several fundamental papers 
he determined the structure of the graded  ring of Siegel modular forms
of degree  two, cf.  \cite{Ig2}, \cite{Ig3} and \cite{Ig4} 
using different methods. 

First he applied his study \cite{Ig1} of moduli theory of 
the (hyperelliptic) curves of genus two in view of binary sextics,
cf. \cite{Ig1}.
The second method is based on, what he calls, a 'fundamental lemma'
\cite{Ig}.
Using successively invariant theory of a finite group 
to the fundamental lemma,
he determined the ring of modular forms of degree two, cf. \cite{Ig3}.
In the third method \cite{Ig4},
he defines the $\rho$ homomorphism from a subring of the ring of 
modular forms to the projective invariants of binary forms.
Then classical invariant theory is used to obtain the ring of
modular forms.


Few years later in \cite{Fr},  Freitag used a more analytic  method
(see also \cite{Hm}): 
he
did a careful analysis of the  
vanishing locus of a distinguished modular form. 
This allowed him to reduce the problem to the study of graded ring
related to lower dimensional varieties.\smallskip

After this,
only a few other cases were determined in the case of Siegel modular
forms of degree two, because of the difficulty of the subject, cf.
\cite{Ib} and \cite{FSM}.

Tsuyumine, cf. \cite{Ts}, was the first one 
who determined the generators of the graded  ring of 
Siegel modular forms of degree three. 
He got the result passing through a lot of difficulties. 
Essentially he used a mix between Igusa's third method and 
Freitag's method and deeply used Shioda's result 
about the structure of the ring of binary octavics, cf. \cite{Sh}.
\smallskip

Later Runge used Igusa's second method for determining the graded ring, 
in fact  using the result of \cite{SM}, 
he was able to determine the structure of a graded  ring of 
Siegel modular forms of degree  three related to a subgroup of 
finite index  of the integral symplectic group. 
Then he applied Igusa's going down process.\smallskip

Runge's method had also the advantage of relating the ring of 
Siegel modular forms with the ring $R_3$
of code polynomials of genus three. 
In fact he proved that the ring of modular forms is isomorphic to
$R_3/\< J\3 \>$. 
Here $J\3$ is the 
difference of the code polynomials for the $e_8\oplus e_8$ 
and $d_{16}^{+}$ codes.  

Moreover he proved that this map exists for any genus $g$ and 
has the image contained in the graded subring of modular forms of even
weight. Also with this restriction, 
it can be easily shown that this map is not surjective for large $g$. 

But there was always the hope that this map  was surjective 
for other mall values of $g$. In fact, 
since the structure of the ring of code polynomials is easier to
determine  
than the structure of the ring of modular forms, 
one would like to apply Runge's method to higher $g$.\smallskip

However, a previous result of the second named author, cf. \cite{SM2},
implied that the map is not surjective when $g\geq 5$.  
Still it was open in the case $g=4$.
Moreover if one restrict the  attention to code polynomials of degree 
divisible by 8 and hence to Siegel modular forms of 
weight divisible by 4, 
the negative result in \cite{SM2} can be by passed.\smallskip

However in \cite{OPY} was observed that as a consequence of the results 
of \cite{E} or  \cite{BFW} and \cite{NV} 
this map could not be sujective if $g\geq6$. 
In fact  dimension of the space of modular forms of weight 12 is 
greater than the dimension of  space of weight enumerators of degree 24.

Meanwhile in genus 4 there were some partial results in \cite{FO} and in \cite{OPY} leading toward a possible surjectivity.\smallskip

In this  note we  show that the map
$$ Th_2: R_g^{(8)}\to M(\Gamma_g)^{(4)}$$
is not surjective if $g\geq 4$.\smallskip

Here $M(\Gamma_g)$ stands for the graded  ring of Siegel modular forms 
of degree  $g$ and the exponents $(8)$ and  $(4)$ mean  
that we consider only code polynomials of degree divisible by 8 and
modular forms of weight divisible by 4 that are related to theta series.

\section{ Notation and Context}

Let $g$ be a positive integer.
We denote by $\Gamma_g:=Sp(g,\Z)$ the integral symplectic group;
it acts
on the Siegel upper-half space $\H_g$ by
$$
\s\cdot \tau:=(A\tau+B)(C \tau+D)^{-1},
$$
where $\s=\pmatrix A&B\\ C&D\endpmatrix\in\Gamma_g.$ 
An element $\tau$ of $\H_g$
is called {\it reducible} if there exists
$\s\in\Gamma_g$ such that
$$
\s\cdot\tau=\pmatrix \tau_1&0\\
0&\tau_2\endpmatrix,\quad\tau_i\in\H_{g_i},g_i>0,\ \ g_1+g_2=g;
$$
otherwise we say that $\tau$ is irreducible.\medskip

Let $k$ be a positive integer and $\Gamma$ be a finite index subgroup.
A multiplier system of weight $k/2$ for $\Gamma$ is 
a map $v:\Gamma\to \C^*$, such that the
map
$$
  \sigma\mapsto v(\sigma)\det(C\tau+D)^{k/2}
$$
satisfies the cocycle condition for every $\sigma\in\Gamma$ and
$\tau\in\H_g$ (note that the function $\det(C\tau+D)$ possesses a
square root). \smallskip

We shall write
$f\vert_{r/2}\s$ for $\det(C\tau+D)^{-k/2}f(\s\cdot\tau)$.

With these notations,  we say that
a holomorphic function $f$ defined on $\H_g$ is  a modular form 
of weight $k/2$ with respect to $\Gamma$ and $v$ if 
$$
  f\vert_{r/2}\s=v(\s) f \qquad 
  \forall\sigma\in\Gamma,
$$
and if additionally $f$ is holomorphic at all cusps when $g=1$.\smallskip

We denote by $[\Gamma, r/2, v]$ the  vector space of such functions.
We shall consider  the graded ring
$$M (\Gamma, v):=\bigoplus_{k=0}^{\infty}[\Gamma,k, v^k]. $$
We omit the multiplier if it is trivial.
\medskip

For $w\in\C$ write $\mathbf{e}(w)=e^{2\pi i w}$.
For $\tau \in\H_g$ and column vectors
$z\in\C^g,a,b\in\Z^g$,
we define
the theta function
by
$$
\tt a b ( \tau, z ) =
\sum_{m\in \Z^g} \mathbf{e}\left({
\frac12(m+a)'\tau (m+a)+(m+a)'(z+b)  }\right),
$$
where $X'$ denotes the transpose of $X$.

For any $a\in \Z^g$, the function
$\theta_r[a]$ defined by
$$\theta_r[a](\tau )=
\t{\left[\begin{matrix}{a/r}\\ {0}\end{matrix}\right]}(r\tau, 0 )$$
is  called an
$r$-th order theta-constant. They depend only on $a$ modulo $r$.
 \smallskip

For any even positive integer $r$, 
we denote by $\Gamma_g(r,2r)$
the subgroup of $\Gamma_g$ of elements $\sigma$
satisfying
$$\sigma=\pmatrix  A & B \\  C & D \endpmatrix
\equiv 1_{2g}\mod r,$$ 
$$(B)_0\equiv  (C)_0 \equiv 0\mod 2r,
$$
where $(X)_0$ means to take
the vector determined by the diagonal coefficients of
a square matrix $X$.
If we drop the second condition,
we get the principal congruence subgroup $\Gamma_g (r)$ of level $r$.\smallskip

The functions
$\theta_r[a]$ belong to $[\Gamma_g(r,2r), 1/2,  v_r]$ for a suitable multiplier $v_r$.\smallskip

We denote by $\vec\theta_r=\left[\theta_r[a]
\right]_{a\in  \Z^g/ r\Z^g}$
the vector of $r$-th order theta-constants.

We define $r^g$ variables $F_a$ for $a \in \Z^g/r\Z^g$.
Let $\C[F_a: a \in \Z^g/r\Z^g]$ be the polynomial
ring in these variables and $\C[F_a: a \in \Z^g/r\Z^g]^{(2)}$ the
subring of even degree.
For even $r$, there is a  theta map
$$
\Thr: \C[F_a: a \in \Z^g/r\Z^g]^{(2)}  \to M(\Gamma_g(r,2r),\,\chi)
$$
induced by sending  $F$ to $\vec\theta_r$ and
the ring $M(\Gamma_g(r,2r),\,\chi)$ is the integral closure of
$\IM\Thr$ inside its quotient field when $r$ is greater than or 
equal to $4$, cf. \cite{Ig} and \cite{Mu}. 
Since we are considering modular forms of integral weight, 
the multiplier $\chi$ is a character.
It is trivial if and only if $4$ divides $r$, 
otherwise $\chi^2$ is trivial.\smallskip

This theorem is called
the `fundamental lemma' of Igusa.\smallskip

The second named author  proved the same conclusion in the case  $r=2$, cf. \cite{SM} and \cite{Ru}
for the  theta map
$$
\Thtwo: \C[F_a: a \in \F_2^g]^{(2)}  \to M(\Gamma_g (2,4), \chi)
$$
Also in this case $\chi^2$ is trivial. Let $\Gamma_g^* (2,4)$ be the
kernel of  $\chi$. 
This is  obviously
a subgroup of index two in $\Gamma_g (2,4)$, described by the condition
$tr(A-1_g)\equiv 0\pmod{4}$, cf. \cite{Ru}.\smallskip

To the map $\Thtwo$ is associated a projective map
$$ th_2:\overline{\H_g/\Gamma_g (2,4)}\to\P^{2^g -1}.$$
Here with the bar we denote  the Satake compactification of the modular
variety.
This map has been studied in details in \cite{SM}.
In this paper is first proved that the map is generically injective,
then with an ad hoc argument is proved the injectivity.  
This second 
part
appears at the author incomplete, so at the moment we can say
that the map is generically injective. 
However Runge's results in
\cite{Ru} imply that 
when $ g\leq 3$, the map is injective.\smallskip

We know that
the group
$\Gamma_g$ is generated by the elements
$J=\pmatrix 0&-1\\ 1&0\endpmatrix $ and
$t(S)=\pmatrix 1&S\\ 0&1\endpmatrix  $ for integral symmetric $S$.
Moreover it acts on the $r$-th order theta-constants. 
On the second order 
we recall the action of the generators, cf. \cite{Ru} for details.
For these elements, we have
$$
\vec{\theta}_2 \vert_{\frac12}
t(S)=
D_S \vec{\theta}_2
\text{ and }
\vec{\theta}_2 \vert_{\frac12}
J=
\pm\, T_g\vec{\theta}_2,
$$
where
\[
 D_S=\diag (i^{a'Sa})_{a\in \F_2^g}
\text{ and }
T_g=\left(\frac{1+i}{2}\right)^g\left((-1)^{(a,b)}\right)_{a,b\in \F_2^g}.
\]
The $\pm 1$ comes from the choice of square root.
The group $H_g=\< T_g,\{D_S\}_S\>\subseteq \GL(2^g,\C)$    
is of finite order and the representation $\phi:\Gamma_g\to H_g/\pm 1$
defined by $\vec{\theta}_2\vert\s=\pm \phi(\s)\vec{\theta}_2$
defines $\Gamma_g^{*}(2,4)$ as its kernel.  
Hence the map $th_2$ results to be $ \Gamma_g/\Gamma_g(2,4)$ equivariant.
We observe that we consider the $ \Gamma_g(2,4)$ quotient of $
\Gamma_g$, 
since being the map 
projective, twisting by  character we do not produce 
any changement in the image.
\smallskip

We see that 
an $H_g$-invariant polynomial goes to a level
one Siegel form  of even weight under the map $\Thtwo$.
We denote by $R_g$ the $H_g$-invariant subring of
$ \C[F_a: a \in \F_2^g]^{(2)}$ and $R_g^{m}$ the
vector space of  $H_g$-invariant homogeneous polynomials
of degree $m$, thus we have a  theta map
$$\Thtwo:R_g\to M(\Gamma_g)$$
whose image is contained in $M(\Gamma_g)^{(2)}$, 
i.e., in the subring of modular forms of even weight.\smallskip

We want to consider also the action of the group $\Gamma_g$ 
on the fourth order theta-constants.
It is useful to consider this action on different representatives 
that we are going to define.\smallskip

A characteristic $m$ is a column vector in $\Z^{2g}$, 
with $m'$ and $m''$ as first and second entry vectors. 
We put $$e(m)=(-1)^{(m', m'')}$$ 
and we say that $m$
is even or odd according as $e(m)=1$ or $-1$.
Here we denoted by $(\cdot\, ,\cdot)$ the standard scalar product

We  consider the  theta function
$$\vartheta_m( \tau, z)=\t{\left[\begin{matrix}{m'/2}\\ 
{m''/2}\end{matrix}\right]}( \tau, z ).$$ 
This is also called even or odd according if $m$ is  even or odd.
We recall from \cite{Ig0} the following formula:
$$\t{\left[\begin{matrix}{m'/4}\\ {0}\end{matrix}\right]}
(4\tau, 4z )=\frac{1}{2^g}\sum_{m''}\mathbf{e}
((1/2) (m',m'') )\vartheta_m(\tau, 2z).$$
Thus we can  consider the theta-constants 
$\vartheta_m:=\vartheta_m (\tau,0)$ as the entries of the vector 
$\vec\theta_4$. 
We remind that in this case the entries 
at 
odd characteristics are $0$.\smallskip

On these entries the action of $\Gamma_g$ is simpler, in fact it is monomial,
cf. \cite{Ig0} or \cite{Ig2}.\smallskip

Moreover we recall the addition  formula relating theta-constants of 
the second order with theta-constants of the fourth order:

$$\vartheta_m^2=\sum_{a\in \F_2^g} 
(-1)^{(a,m'')}\theta_2[a+m']\theta_2[a].$$

\section{The map $th_2$}

In this section we shall consider  with more details the map $th_2$.
In particular we shall consider its injectivity on some special subloci.  
\smallskip

We recall that in  \cite{Sa} has been proved that the map is injective 
along the hyperelliptic locus. \smallskip

Here we are interested in the locus of the completely
reducible periods, i.e. to the points
$\tau$ that are $\Gamma_g$ conjugated to points of the type
$$\left(\begin{array}{ccccc}
  \lambda_1&0&\dots&\dots&0 \\
0&\lambda_2& 0 &\dots & 0\\ \dots&\dots&\dots&\dots&\dots\\
\dots&\dots&\dots&\dots&\dots\\
0 &\dots &\dots &0&\lambda_g\\
    \end{array}\right) \qquad \lambda_i\in \H_1
  .$$
For doing this we need to  recall and improve some results 
in \cite{SM}.\smallskip

To any characteristic $m$ is associated a character $\chi_m$ of the group
 $G'=\Gamma_g (2,4)/\Gamma_g(4,8)$. We set $G=G'/\pm 1$, then the group of characters  $\hat G$ is  spanned by the characters $\chi_m$ associated to the even characteristics.\smallskip

For any point $\tau_0\in\H_g/\Gamma_g(4,8)$ 
we consider the subgroup  $H_{\tau_0}$ of $\hat G$ spanned by 
all characters $\chi_m\chi_n$ such that 
the products $\vartheta_m\vartheta_n$ do not
vanish at $\tau_0$. Moreover we set $St_{\tau_0}$ be 
the subgroup of $G'$ generated by all $\s$ fixing the point $\tau_0$.

We recall from \cite{SM} the  following

 \begin{thm} The following statements are equivalent:

 (i) $\tau_0$ is a reducible point.

 (ii) $St_{\tau_0}$  is different from $\pm 1$.

 (iii)  $H_{\tau_0}$  is a proper subgroup of  $\hat G$.

 \end{thm}

We need to mention some facts about this theorem.  
In the proof the author uses the theorem stating that 
the map $\th_2$ is injective. This, in general,  is not necessary 
in fact theta-constants are used only to separate points in the 
same fibre  of the  covering map
$$\pi:\H_g/\Gamma_g(4,8)\to \H_g/\Gamma_g(2,4).$$ 

Moreover, since $St_{\tau_0}/\pm 1$ and $H_{\tau_0}$ are dual 
we have that the number of   points in the fiber of $\pi(\tau_0)$ 
is equal to the order of $H_{\tau_0}$.\smallskip

Finally, cf.\cite{SM4} also, we have the following

 \begin{cor} Let $\tau_0\in \H_g/\Gamma_g(4,8)$. 
It is $\Gamma_g$ equivalent to a point of the form
 $\tau$ that are $\Gamma_g$ conjugated to points of the type
$$\left(\begin{array}{ccccc}
  \lambda_{i_1}&0&\dots&\dots&0 \\
0&\lambda_{i_2}& 0 &\dots & 0\\ \dots&\dots&\dots&\dots&\dots\\
\dots&\dots&\dots&\dots&\dots\\
0 &\dots &\dots &0&\lambda_{i_k}\\
    \end{array}\right) 
  $$
with all $\lambda_{i_j}\in\H_{g_j}$ irreducible and $g_1+g_2+\dots + g_k=g$
if and only if  $St_{\tau_0}$  has order $2^k$.
 \end{cor}

Let $\tau_0$ be a  completely reducible point, 
already in diagonal form.
It is immediate to verify that
 $St_{\tau_0}$ is the subgroup of $G'$  that is the image of $\s\in\Gamma_g(2,4)$ that are diagonal matrices, i.e. $B=C=0$, $A=D= diag(\pm1,\pm1,\dots,\pm1)$. Obviously it has order $2^g$.

  \begin{prop} 
The map $th_2$ is injective along the completely reducible points.
 \end{prop}

 {\it Proof}. 
Assume that 
 $$  th_2(\tau)= th_2(\tau')$$
for a completely reducible point $\tau$ and another point $\tau'$.
Here we can take $\tau$ in diagonal form.
By addition formula,
the same first order theta-constants 
vanish on $\tau$ and on $\tau'$.
Thus $\tau'$ is completely reducible and 
conjugate to $\tau$ by an element $\s$ in the group $\Gamma_g(2) $ 
that fixes the characteristics of the first order theta-constants. 
But now $\Gamma_g(2) $ acts inducing  different characters 
on the $\theta_m^2$. 
Requiring that the character has to be the same, 
we have that $\s\in \Gamma_g(2, 4)$.

Another interesting fact about the completely reducible points is 
the following 

\begin{prop} 
Let $\tau_0$ be a generic completely reducible point.
If $\s\in \Gamma_g$ stabilizes  $\tau$, i.e.
  $$\s\tau_0=\tau_0,$$
then $\s\in \Gamma_g(2,4)$.
  \end{prop}

{\it Proof}. 
Without any loss of generality we can assume $\tau_0$ in diagonal form with
$\tau_0= x_0+iy_0$, with $x_0$ and $y_0>0$ diagonal. 
If $\s$ stabilizes $\tau_0$ we have    
   $$A\tau_0 +B= \tau_0(C\tau_0 +D).$$
>From this we get
$$ Ax_0+B=x_0(Cx_0 +D) -y_0Cy_0 $$
and 
  $$Ay_0= y_0(Cx_0+D) + x_0Cy_0.$$
  We can choose $x_0$ and $y_0$ in such form  that  for any $\s\in \Gamma_g$
  $$Ax_0-x_0(Cx_0 +D) +y_0Cy_0$$ is not an integral matrix, 
hence necessarily we have $B=0$.
Moreover from the first equation we can choose $y_0>>0$, 
so that  $C=0$, hence we have
   $$ Ax_0 =x_0D \quad {\rm and}\quad   Ay_0= y_0D .$$
These conditions for generic diagonal matrices imply $A=D$ diagonal  and
  integral. 
Thus we get the desired result. \medskip

We conclude this section describing non-embeddability 
of $\H_g/\Gamma_g(2,4)$.

First we recall that in the paper
\cite{SM3} has been proved the following result.

\begin{prop}  
When $g\geq 4$, the map $th_2$  is not an immersion 
at the completely reducible points.
\end{prop}

In fact  the second author computed explicitly the dimension 
$t_g$ of the tangent spaces 
at the generic points completely reducible points of $\H_g/\Gamma_g(2,4)$.
Since there is a misprint in the formula,
we reproduce it here:
$$t_g=g+ {g\choose 2} +\frac{1}{2}\sum_{h=3}^g {g\choose h} (h-1)!,$$
in which the third term is read as zero when $g=1,2$.
  
We recall the idea of the proof of this formula.
Let $K $ be the subgroup of $GL(g, \Z)$ formed 
by the diagonal matrices and $S_g(\C)$ be the set of symmetric matrices, 
thus a neighborhood of a generic completely reducible point 
$\tau_0\in\H_g/\Gamma_g(2,4)$ looks like a neighborhood of 
$0\in S_g(\C)/K$.\smallskip

Moreover the ideal ${\mathfrak m}_0 \subset \C[X_{i,j}]^K$ is 
generated by the monomials
  
$$X_{11},X_{22},\ldots,X_{gg},\text{ and }
X_{i_1, i_2}X_{i_2, i_3}\dots X_{i_n, i_1}$$
with $1\leq i_1<i_2<\dots<i_n\leq g$.\smallskip
 
These monomials are also a basis of ${\mathfrak m}_0/{\mathfrak m}_0^2$. 
Thus the dimension of the tangent space is the  number of 
such monomials that is exactly $t_g$.
 
When $g\geq 4$, we have $t_g>2^{g}-1$ and hence 
we cannot have an immersion of  a neighborhood of the point $\tau_0$
into $\P^{2^g -1}$.
We have thus proved 
that $\H_g/\Gamma_g(2,4)$ cannot be embedded in $\P^{2^g -1}$
when $g\geq 4$.

\section{ The map $\Theta_2$}

In this section we prove our main result.
As we mentioned before,
the transformation formula of theta-constants implies 
that the map $\th_2$ is $\Gamma_g/\Gamma_g(2,4)$ equivariant.
Therefore we have a map
 $$\Thtwo: R_g\to M(\Gamma_g)$$
 whose image is contained in $ M(\Gamma_g)^{(2)}$.\smallskip

>From Runge's results, cf. \cite{Ru}, \cite{Ru2} and \cite{Ru3}, 
we have that the map is surjective onto $ M(\Gamma_g)^{(2)}$ 
when $g\leq 3$.\smallskip

In general
a map $\Thtwo$ induces a projective map
 $$\Theta_2: \overline{\H_g/\Gamma_g}\to Proj( R_g).$$
 As a consequence of Runge's results it is an immersion when $g\leq 3$.
 For higher $g$, we have the following

\begin{thm} 
$\Theta_2: \overline{\H_g/\Gamma_g}\to Proj( R_g) $
is not an embedding when $g\geq 4$
\end{thm}
{\it Proof}.
We consider the action of $K_g=\Gamma_g/\Gamma_g(2,4)$ 
at a generic completely reducible point $\tau \in \H_g/\Gamma_g(2,4)$. 
The result of Proposition $4$ implies that the group acts freely on
$\tau$. 
Let $x=\th_2(\tau)$, since the map is equivariant the group acts 
freely also on $\P^{2^g-1}$.\smallskip

Let
$\pi:  \H_g/\Gamma_g(2,4)\to  \H_g/\Gamma_g$ and 
$\phi: \P^{2^g-1}\to \P^{2^g-1}/K_g$, 
the maps that make commutative the following diagram

  \xymatrix{
    &  \H_g/\Gamma_g(2,4)\ar[d]^{\pi}\ar[r]  &\P^{2^g-1}\ar[d]^{\phi}\\
&\H_g/\Gamma_g \ar[r] & \P^{2^g-1}/K_g.}
We set $\tau ':=\pi(\tau)$ and $x'=\phi(x)$. 
Since the action of the group was free we have that 
the dimension of the tangent spaces at $\tau '$ and at $x'$ are equal 
to the dimension of the tangent
spaces at $\tau$ and $x$, respectively. 
Hence, in particular, they have different  dimension, thus 
$\Theta_2$ is not an immersion. \bigskip

As an immediate consequence, using basic fact from \cite{Ha}, page 92 , Exercise 3.12,  we have that

\begin{cor} When $g\geq 4$,  for any $k$ the homomorphism
   $$\Thtwo: R_g^{(4k)}\to M(\Gamma_g)^{(2k)}$$
is never surjective.
\end{cor}

We interpret this corollary in a restricted case,
which is of another interest.
To do this,
we denote by $\C [\vartheta_{\Lambda}]$ the ring of theta series
of all even unimodular lattices $\Lambda$ 
where
\[
 \vartheta_{\Lambda}(\tau)
=\sum_{v_1,v_2,\ldots,v_g\in \Lambda}
\prod_{i,j} \mathbf{e} ((v_i,v_j)\tau_{ij}/2).
\]
Note that an even unimodular lattice exists if and only if
the rank of a lattice is a multiple of $8$
and that the weight of the corresponding theta series
is a half of the rank.
We recall that for $ n$ large enough, 
cf. \cite{Bo} or \cite{We}, 
the space of modular forms $[\Gamma_g,n]$ is spanned by theta series 
provided that $4$ divides $n$.
Our final result in this paper is, 
specializing the previous corollary to the case $k=2$, 

\begin{cor} 
When $g\geq 4$ the homomorphism
$$\Thtwo: R_g^{(8)}\to\C[\vartheta_{\Lambda}]$$
is not surjective for infinitely many degrees.
\end{cor}




\begin{thebibliography}{99}

\bibitem{BFW}
Borcherds, R, E., Freitag, E., Weissauer, R.,
A Siegel cusp form of degree 12 and weight 12,
J. Reine Angew. Math. 494 (1998), 141--153.




\bibitem{Bo}
B\"{o}cherer, S.,
\"{U}ber die Fourier-Jacobi-Entwicklung Siegelscher Eisensteinreihen,
Math. Z. 183 (1983), no. 1, 21--46.




\bibitem{E}
Erohin, V. A.,
Theta series of even unimodular $24$-dimensional lattices,
Algebraic numbers and finite groups.
Zap. Nau\v cn. Sem. Leningrad. Otdel. Mat. Inst. Steklov. (LOMI) 86
(1979), 82--93, 190.


\bibitem{Fr}
Freitag, E.,
Zur Theorie der Modulformen zweiten Grades,
 Nachr. Akad. Wiss. G\"ottingen, II. Math.-Phys. Kl. 1965, 151-157 (1965).

\bibitem{FO}
Freitag, E, Oura, M.,
A theta relation in genus 4,
Nagoya Math. J. 161 (2001), 69--83.


\bibitem{FSM}
Freitag, E., Salvati Manni, R.,
The Burkhardt group and modular forms,
Transform. Groups 9 (2004), 25--45.

\bibitem{Ha}
Hartshorne, R.
Algebraic Geometry,
GTM 52, Springer- Verlag, New York Heidelberg Berlin  1977.




\bibitem{Hm}
Hammond, W. F.,
On the graded ring of Siegel modular forms of genus two,
Amer. J. Math. 87 1965 502--506.


\bibitem{Ib}
Ibukiyama, T.,
On Siegel modular varieties of level $3$,
Internat. J. Math. 2 (1991), no. 1, 17--35.

\bibitem{Ig0}
Igusa, J.,
Theta functions,
Die Grundlehren der mathematischen Wissenschaften, Band 194. 
Springer-Verlag, New York-Heidelberg,1972.


\bibitem{Ig1}
Igusa, J.,
Arithmetic variety of moduli for genus two,
Ann. of Math. (2) 72 1960 612--649.




\bibitem{Ig2}
Igusa, J.,
On Siegel modular forms of genus two,
Am. J. Math. 84(1962), 175-200.


\bibitem{Ig3}
Igusa, J.,
On Siegel modular forms of genus two (II),
Am. J. Math. 86(1964), 392-412.


\bibitem{Ig}
Igusa, J.,
On the graded ring of theta-constants,
Amer. J. Math. 86 (1964), 219--246.

\bibitem{Ig4}
Igusa, J.,
Modular forms and projective invariants,
Amer. J. Math. 89 (1967), 817--855.

\bibitem{Mu}
Mumford, D.,
Tata lectures on theta III,
with collaboration of M. Nori and P. Norman, 
Birkh\"{e}user Boston, Inc.,
Boston, MA, 2007.



\bibitem{NV}
Nebe, G., Venkov, B.,
On Siegel modular forms of weight 12,
J. Reine Angew. Math. 531 (2001), 49--60.

\bibitem{OPY}
Oura, M., Poor, C., Yuen, D.S.,
Toward the Siegel ring in genus four,
to appear in Int. J. Number Theory.

\bibitem{Ru}
Runge, B.,
On Siegel modular forms. I,
J. Reine Angew. Math. {\bf 436}, 57-85 (1993).


\bibitem{Ru2}
Runge, B.,
On Siegel modular forms II,
Nagoya Math. J. 138 (1995), 179--197.

\bibitem{Ru3}
Runge, B.,
Codes and Siegel modular forms,
Discrete Math. 148 (1996), no. 1-3, 175--204.



\bibitem{SM2}
Salvati Manni, R.,
On the not integrally closed subrings of the ring of the thetanullwerte,
Duke Math. J. 52 (1985), no. 1, 25--33.


\bibitem{SM}
Salvati Manni, R.,
Modular varieties with level $2$ theta structure,
Amer. J. Math. 116 (1994), no. 6, 1489--1511.


\bibitem{SM3}
Salvati Manni, R.,
On the projective varieties associated with some subrings of the ring of 
Thetanullwerte,
Nagoya Math. J. 133 (1994), 71--83.



\bibitem{SM4}
Salvati Manni, R.,
On the locus of reducible abelian varieties,
Workshop on Abelian Varieties and Theta Functions (Morelia, 1996), 151--155,
Aportaciones Mat. Investig., 13, Soc. Mat. Mexicana, 1998.

\bibitem{Sa}
Sasaki, R.,
Modular forms vanishing at the reducible points of the Siegel upper-half space,
J. Reine Angew. Math. 345 (1983), 111--121.


\bibitem{Sh}
Shioda, T.,
On the graded ring of invariants of binary octavics,
Amer. J. Math. 89 (1967), 1022--1046.

\bibitem{Ts}
Tsuyumine, S.,
On Siegel modular forms of degree three,
Am. J. Math. 108 (1986), 755-862; Addendum 1001-1003 (1986).


\bibitem{We} 
Weissauer, R.,
Stabile Modulformen und Eisensteinreihen,
Lecture Notes in Mathematics, 1219. Springer-Verlag, Berlin, 1986.


\end{thebibliography}
 \end{document}